\documentclass[12pt]{article}

\usepackage[margin=1in]{geometry}
\usepackage{amsmath, amssymb, amsfonts}
\usepackage{amsthm}
\usepackage{mathtools}
\usepackage{bm}
\usepackage{cite}
\usepackage[colorlinks=true,linkcolor=blue,citecolor=blue,urlcolor=blue]{hyperref}

\numberwithin{equation}{section}
\newtheorem{theorem}{Theorem}[section]

\theoremstyle{definition}
\newtheorem{definition}[theorem]{Definition}
\theoremstyle{remark}

\newcommand{\N}{\mathbb{N}}
\newcommand{\Z}{\mathbb{Z}}
\newcommand{\R}{\mathbb{R}}
\newcommand{\C}{\mathbb{C}}

\renewcommand{\P}{\mathbb{P}}
\newcommand{\E}{\mathbb{E}}

\newcommand{\Zf}{\mathbb{Z}_{\mathrm{f}}}

\newcommand{\Cyln}{\mathcal{C}_n}     

\title{Classical localization problem: a survey}
\author{Zoey Zhou
\thanks{Carnegie Mellon University, Pittsburgh, PA.} }
\date{}

\begin{document}
\maketitle

\begin{abstract}
We survey classical localization problems arising from quantum network models in symmetry class~C and their mappings to history–dependent random walks on directed lattices. We describe how localization versus delocalization of trajectories can be analysed using percolation methods and combinatorial enumeration of path intersection patterns. In particular, we review results establishing almost sure finiteness of trajectories for parameters near criticality and polynomial bounds on the confinement length in cylindrical geometries.
\end{abstract}

\section{Introduction}

Disordered systems exhibit a striking trichotomy of large--scale dynamical behaviours: ballistic transport, localization, and diffusive transport. In the quantum setting, the paradigmatic example is the Anderson tight--binding model on $\Z^d$, introduced in Anderson’s seminal paper~\cite{And58}, in which a discrete Laplacian is coupled to an i.i.d.\ random potential and one observes the \emph{absence of diffusion in certain random lattices}. Since Anderson’s work, a vast body of physics and mathematics has developed around quantum localization, scaling theory, and the metal--insulator transition; see, for instance, the original scaling theory paper of Abrahams--Anderson--Licciardello--Ramakrishnan~\cite{AALR79}, the review of disordered electronic systems by Lee--Ramakrishnan~\cite{LeeRamakrishnan1985}, and the modern overview in Evers--Mirlin~\cite{EversMirlin2008}. On the mathematical side, there is now a detailed theory of Anderson localization for a variety of single--site distributions, including singular ones; see, for example \cite{CKM87,BK05,LiZhang3D,Li20AB2D}.

From the physics perspective, a central theme is \emph{universality}: the large--scale behaviour near the localization--delocalization transition should depend only on the spatial dimension and the symmetry class of the model, but not on microscopic details. This unifying principle appears across statistical mechanics and random media. In classical percolation, macroscopic connectivity undergoes a sharp phase transition, and critical exponents are expected to be universal within broad classes of planar graphs; see Grimmett’s monograph~\cite{Grimmett99}. In quantum disordered systems, nonlinear sigma models, the supersymmetry method, and random matrix theory underpin the conjectural universality classes of Anderson transitions~\cite{Wegner1979,Hikami1981,Zirnbauer1996,Efetov1983,EfetovBook,AltlandZirnbauer1997,EversMirlin2008}. In particular, the field--theoretic description of localization in terms of nonlinear sigma models, initiated by Wegner and Hikami~\cite{Wegner1979,Hikami1981} and put on a supersymmetric footing by Efetov~\cite{Efetov1983,EfetovBook}, provides a unifying framework in which symmetry classes appear as target spaces of the sigma model, and different universality classes of Anderson transitions can be classified via Cartan’s list of symmetric spaces~\cite{Zirnbauer1996,AltlandZirnbauer1997}.

A key idea, heavily exploited in the physics literature, is that suitably averaged quantum transport in certain symmetry classes can be rephrased in terms of \emph{classical} stochastic geometry, notably percolation and random line models. The most famous example of this philosophy is the Chalker--Coddington network model for the integer quantum Hall effect~\cite{KramerOhtsukiKettemann2005,Ludwig1994,Huckestein1995}, where critical quantum transport at plateau transitions is mapped to the geometry of classical directed paths undergoing scattering at saddle points. Closely related constructions appear in the sigma--model description of quantum Hall and spin quantum Hall transitions~\cite{Hikami1981,Khmelnitskii1983,Huckestein1995,EversMirlin2008}.

\medskip

A particularly clean realization of such a quantum--to--classical reduction is due to Cardy in the context of \emph{class C} network models, which arise from quasiparticle dynamics in disordered spin--singlet superconductors with broken time--reversal symmetry but preserved spin--rotation invariance. On the one hand, one considers a unitary evolution on a directed graph $G=(V,E)$: each edge $e\in E$ carries a complex amplitude, and at each vertex $v\in V$ a random $2\times 2$ scattering matrix redistributes incoming amplitudes among outgoing edges, subject to the class~C symmetry constraints determined by Altland--Zirnbauer’s classification~\cite{AltlandZirnbauer1997}. The randomness is implemented by choosing these unitary scattering matrices independently from an ensemble invariant under the relevant symmetry group, in close analogy with the random matrix ensembles used in nuclear physics and mesoscopic systems~\cite{AltlandZirnbauer1997,EversMirlin2008}. On the other hand, one can introduce a \emph{classical} model on the same graph: configurations consist of collections of disjoint, directed paths and loops following the edge orientation, subject to local pairing rules at each vertex.

Cardy’s argument shows that, after averaging over the quantum disorder using supersymmetric techniques in the spirit of Efetov~\cite{Efetov1983,EfetovBook}, certain Green–function observables of the network model are exactly equal to correlation functions of the classical line model~\cite{CardyClassC,BOC03,Car10}. More concretely, if $G$ is a planar, oriented graph and $\mathcal U$ denotes the random one–step unitary evolution, then for appropriate observables one has identities of the schematic form
\[
  \mathbb E\big[\,|\langle \delta_x, \mathcal U^t \delta_y\rangle|^2\,\big]
  \;=\; \P_{\mathrm{cl}}\big(\text{there is a classical path from $x$ to $y$ of length $t$}\big),
\]
where the right–hand side is taken with respect to a probability measure on classical paths and loops determined by the same microscopic disorder~\cite{CardyClassC,BOC03}. In particular, the long–time behaviour of the quantum evolution---such as localization versus delocalization of wave packets---is reduced to a question about the geometry of (random) classical trajectories. The scaling theory picture, originating in~\cite{AALR79} and refined in the context of quantum Hall transitions~\cite{Khmelnitskii1983,Ludwig1994,Huckestein1995}, predicts that in two dimensions all quantum states are localized in the corresponding random Schr\"odinger model in the absence of special topological mechanisms; on the classical side this corresponds to the conjecture that in certain planar network models \emph{all} classical trajectories are almost surely finite for every positive disorder parameter~\cite{BOC03,Spencer2012}.

\medskip

On the Manhattan lattice, this reduction becomes especially transparent. The Manhattan lattice is the oriented graph on $\Z^2$ obtained by orienting horizontal lines alternately east and west and orienting vertical lines alternately north and south. Beamond, Owczarek and Cardy~\cite{BOC03} considered a class~C network model on this lattice and showed that, after disorder averaging, one is led to a classical model for trajectories on a percolation configuration: the \emph{Manhattan pinball problem}. In the formulation used in~\cite{Spencer2012,Li21Manhattan}, let $\Zf^2$ denote the tilted square lattice and consider Bernoulli bond percolation with parameter $p\in(0,1)$ on $\Zf^2$. Each closed edge carries a mirror that deflects a light ray through a right angle, while open edges are traversed straight. A ray of light travels along the Manhattan orientation and turns only when encountering a mirror. The central question is: for which values of $p$ is the resulting trajectory almost surely finite? This classical problem is intimately related, via Cardy’s mapping, to localization properties of the underlying quantum network model~\cite{BOC03,CardyClassC,Car10} and hence to the two–dimensional Anderson localization scenario~\cite{AALR79,LeeRamakrishnan1985,EversMirlin2008}.

\medskip

The \emph{Lorentz mirror model} is a closely related classical dynamical system designed as a discrete analogue of the Lorentz gas. Here one considers the nearest–neighbour lattice $\Z^d$, and places two–sided mirrors at vertices with probability $p$, independently. In two dimensions these mirrors are rotated by $\pm \pi/4$ and reflect a light ray deterministically, while in higher dimensions one can formulate the reflection rule in terms of random matchings on the set of coordinate directions~\cite{RuijgrokCohen1988,BOC03}. Launching a ray of light from the origin along an edge, one obtains a deterministic trajectory determined by the mirror configuration. The basic dichotomy is again localization versus delocalization: is the trajectory almost surely finite, or can it travel to arbitrarily large distances? Numerical studies suggest that in $d=2$ all trajectories are finite for every $p>0$~\cite{BOC03,ZiffKongCohen91}, while on certain quasi–one–dimensional geometries extended behaviour may occur~\cite{KS15}. Rigorous progress in $d=2$ has so far been limited to lower bounds on the escape probability~\cite{KS15} and on the density of closed paths~\cite{KraemerSanders14}, whereas more is known on cylindrical geometries~\cite{Li20Cylinder}.

\medskip

The goal of this survey is to present the current state of the rigorous theory for these \emph{classical localization models}, with an emphasis on the technical mechanisms underpinning localization and confined transport. In particular, we focus on two closely related families of models and results:
\begin{itemize}
  \item the Cardy reduction from quantum class~C network models to classical line models on planar graphs, and its implications for universality and the interpretation of classical localization as a proxy for quantum localization~\cite{CardyClassC,BOC03,Car10};
  \item localization on the two–dimensional Manhattan lattice and on a one–dimensional cylinder for the Manhattan pinball and Lorentz mirror models, where percolation methods and combinatorial enumeration play a central role~\cite{Li21Manhattan,Li20Cylinder,Grimmett99,AG91,BalisterBollobasRiordan14,GrimmettSelectedProblems,Zhou25Lattice}.
\end{itemize}

These results rely on a diverse but interconnected toolbox. In the two–dimensional Manhattan pinball problem, the work~\cite{Li21Manhattan} sharpens earlier percolation–based arguments~\cite{BOC03,Spencer2012} and shows that for $p>1/2-\varepsilon$ the trajectory is almost surely bounded. The proof is based on a delicate \emph{enhancement} construction that increases the density of closed edges while preserving percolative structure. This places the problem within the general theory of (non–monotone) enhancement of percolation models developed by Aizenman--Grimmett and Grimmett~\cite{AG91,Grimmett99}, and uses Russo--Seymour--Welsh crossing estimates and Peierls arguments in a crucial way. On the cylinder $\Cyln = \Z \times (\Z/2n\Z)$, the Lorentz mirror and Manhattan models are studied in~\cite{Li20Cylinder}. There, one proves that, for each fixed $p\in(0,1)$, with high probability every trajectory intersecting the section $x=0$ is contained in a strip $|x|\le C n^{10}$. Besides percolation–type control of crossings in tilted rectangles, an essential ingredient is a \emph{combinatorial enumeration} argument that bounds the number of long, left–right trajectories on the cylinder by counting how often they must intersect a set of shorter, up–down closed trajectories and using symmetry to reduce the probability of such configurations.

\medskip

The present article is organized as follows. In Section~\ref{sec:classC} we review quantum network models in symmetry class~C and Cardy’s mapping to classical localization models on planar graphs. We give a precise formulation of the network model on a directed graph, summarize the disorder averaging that leads to a classical measure on random paths, and explain in some detail how the Manhattan lattice realization in~\cite{BOC03,CardyClassC,Car10} reduces questions about quantum localization to the geometry of classical trajectories. We also briefly recall the connection with the two–dimensional Anderson model and scaling–theory predictions~\cite{And58,AALR79,LeeRamakrishnan1985,EversMirlin2008}.

Section~\ref{sec:2d-classical} is devoted to mathematically rigorous localization results in two dimensions and on cylindrical geometries. We first introduce rigorous definitions of the Manhattan pinball and Lorentz mirror models on $\Z^2$ and on the even–width cylinder $\Cyln$, and then survey the percolation–based proof of localization for $p>1/2-\varepsilon$ on the Manhattan lattice~\cite{Li21Manhattan}, emphasizing the enhancement construction, Russo--Seymour--Welsh crossing events, and Peierls–type contour estimates~\cite{AG91,Grimmett99}. We then turn to the polynomial confinement result on the cylinder~\cite{Li20Cylinder}, where we describe how percolation arguments on annuli combine with a combinatorial enumeration of left–right crossings to show that typical trajectories are trapped in a strip of width $O(n^{10})$. Along the way we place these results in the broader context of earlier work on the Lorentz mirror model, including the lower bound on escape probability~\cite{KS15} and numerical studies of open–path densities~\cite{ZiffKongCohen91,KraemerSanders14}.

\section{Quantum network models of class C and Cardy’s mapping}
\label{sec:classC}

In this section we recall the general framework of quantum network models in
symmetry class~C and review Cardy’s mapping from such models to classical,
history–dependent random walks on the underlying graph~\cite{CardyClassC,BOC03,Car10}. We also discuss how this construction fits into the broader physics picture of Anderson localization, sigma models, and symmetry classes~\cite{Efetov1983,EfetovBook,Wegner1979,Hikami1981,Zirnbauer1996,AltlandZirnbauer1997,LeeRamakrishnan1985,EversMirlin2008}. We then specialise to the case of a degree–$4$
planar graph, in particular the Manhattan lattice of~\cite{BOC03}, and explain how questions about quantum localization are reduced to questions about the geometry of classical trails.

\subsection{Network models on directed graphs in class C}

Let $G=(V,E)$ be a finite directed graph. We write $E=\{e_1,\dots,e_M\}$ for the
set of directed edges and $V$ for the set of vertices (nodes). We assume:
\begin{itemize}
  \item each edge $e\in E$ is oriented and connects an ordered pair of vertices;
  \item at each vertex $v\in V$ the number of incoming edges equals the number of outgoing edges; we denote this common number by $N_v$;
  \item the graph is \emph{closed}, in the sense that every edge has both an initial and a terminal vertex in $V$.
\end{itemize}

The one–particle Hilbert space of the network model is
\[
  \mathcal H \;=\; \ell^2(E)\otimes \C^2,
\]
where the factor $\C^2$ carries an internal ``spin'' degree of freedom. This corresponds to the minimal nontrivial representation of the class~C symmetry, in which the Bogoliubov–de~Gennes Hamiltonian has an antiunitary particle–hole symmetry squaring to $-1$ and the effective spin is conserved~\cite{AltlandZirnbauer1997,EversMirlin2008}. We denote by $\{|e,\alpha\rangle: e\in E,\;\alpha\in\{1,2\}\}$ the canonical orthonormal basis of $\mathcal H$, where $e$ labels an edge and $\alpha$ the spin index.

Propagation in the network takes place in discrete time steps. A single time step consists of two parts:
\begin{enumerate}
  \item \emph{Edge evolution}: while travelling along an edge $e\in E$, the wave
  function is multiplied by a random $\mathrm{SU}(2)$ matrix $U_e$ acting on the spin index.

  \item \emph{Node scattering}: when arriving at a node $v\in V$, the particle
  is scattered from incoming edges to outgoing edges according to a fixed real orthogonal matrix $S^{(v)}\in \mathrm{O}(N_v)$.
\end{enumerate}
More concretely, let
\[
  U_{\mathrm{edge}}
  \;=\; \bigoplus_{e\in E} U_e
  \quad\text{on}\quad \mathcal H,
\]
acting trivially on the edge index and as $U_e\in\mathrm{SU}(2)$ on the spin index, and let $U_{\mathrm{node}}$ be the unitary operator which, for each vertex $v$ with incoming edge set $\mathrm{In}(v)$ and outgoing edge set
$\mathrm{Out}(v)$, performs the orthogonal mixing given by $S^{(v)}$ between $\ell^2(\mathrm{In}(v))$ and $\ell^2(\mathrm{Out}(v))$ (and acts trivially on the spin index). The one–step unitary evolution of the network is then
\[
  \mathcal U \;=\; U_{\mathrm{node}}\, U_{\mathrm{edge}}.
\]

The random variables $(U_e)_{e\in E}$ are taken i.i.d.\ with common law the Haar measure on $\mathrm{SU}(2)$. We denote expectation with respect to this disorder by $\E[\cdot]$ or an overline. The node matrices $S^{(v)}\in\mathrm{O}(N_v)$ are kept deterministic; their choice encodes the underlying graph structure and microscopic scattering at the nodes.

Of primary interest is the resolvent (or Green’s function) of $\mathcal U$:
for $z\in\C$ with $|z|\neq1$, define
\begin{equation}
  G(z) \;=\; (1 - z \mathcal U)^{-1}.
  \label{eq:resolvent-def}
\end{equation}
When $|z|<1$ we may expand \eqref{eq:resolvent-def} into a convergent Neumann
series
\[
  G(z) \;=\; \sum_{n\ge0} z^n \mathcal U^n.
\]
The matrix elements of $G(z)$ in the edge–spin basis encode both spectral and
transport properties. Writing
\[
  G_{e,e'}^{\alpha\beta}(z)
  \;=\;
  \langle e,\alpha \mid G(z) \mid e',\beta\rangle,
\]
the \emph{averaged local Green’s function} at edge $e$ is
\[
  \overline{g_e(z)} \;=\; \overline{\operatorname{Tr}_{\C^2}
  G_{e,e}^{\alpha\alpha}(z)},
\]
and the averaged density of states and conductance can be reconstructed from these and related quantities~\cite{Efetov1983,EfetovBook,CardyClassC}. For instance, the averaged density of states for the one–step unitary is
\[
  \rho(\theta)
  \;=\;
  -\frac{1}{2\pi}\,\lim_{\varepsilon\downarrow0}\,
  \Im\,\overline{\operatorname{Tr}_{\mathcal H}
  G(e^{i\theta+\varepsilon})}.
\]

The class~C symmetry implies that the single–particle Hilbert space can be chosen to have even dimension and that the edge evolution is given by $\mathrm{SU}(2)$ matrices, so that after disorder averaging only spin–singlet combinations survive~\cite{CardyClassC,BOC03,EversMirlin2008}. This finite–dimensionality is crucial in what follows: it allows one to express disorder–averaged Green’s functions as sums over \emph{classical} paths with weights depending only on the node matrices $S^{(v)}$.

From the broader localization perspective, this network formulation is a discrete analogue of the nonlinear sigma–model representation of disordered systems~\cite{Wegner1979,Hikami1981,Zirnbauer1996,AltlandZirnbauer1997,Efetov1983}. There, disorder–averaged Green’s functions are expressed as functional integrals over matrix–valued fields taking values in symmetric spaces determined by the symmetry class; in the present context the space of edge variables is finite–dimensional and the integrals factor into ordinary group integrals over $\mathrm{SU}(2)$, which can be evaluated exactly by supersymmetric methods~\cite{Efetov1983,EfetovBook,CardyClassC}. The payoff is that the resulting effective theory is not a continuum field theory but a purely combinatorial measure on trails of a classical walker.

\subsection{Trails and the classical path expansion}

We now briefly describe Cardy’s mapping~\cite{CardyClassC} between the quantum network model and a classical history–dependent random walk on the same graph.
The key step is to express the resolvent \eqref{eq:resolvent-def} as a supersymmetric path integral and to perform the average over the random $\mathrm{SU}(2)$ edge variables $(U_e)_{e\in E}$. This produces a finite expansion in terms of edge–occupancy patterns and, after integrating out all degrees of freedom except a single fermion–fermion (or boson–fermion) pair, leads to a combinatorial representation in terms of certain graph decompositions into loops and open paths.

We only record the outcome, following~\cite{CardyClassC,BOC03}. Fix a directed graph
$G=(V,E)$ as above. At each vertex $v$, label the $N_v$ incoming edges by
$\{1,\dots,N_v\}$ and the $N_v$ outgoing edges by $\{1,\dots,N_v\}$. The node
matrix $S^{(v)}\in\mathrm{O}(N_v)$ has entries $S^{(v)}_{ij}$, where $i$
indexes outgoing and $j$ indexes incoming channels. For a subset
$I\subseteq\{1,\dots,N_v\}$ of incoming channels and a subset
$O\subseteq\{1,\dots,N_v\}$ of outgoing channels with $|I|=|O|=k$, we denote by
\[
  S^{(v)}[O,I]
\]
the $k\times k$ submatrix obtained by restricting $S^{(v)}$ to rows in $O$ and
columns in $I$, in the natural ordered bases, and by
$\det S^{(v)}[O,I]$ its determinant. We write $\mathrm{sgn}(\pi)$ for the
signature of a permutation $\pi$.

\begin{definition}[Trail]
A \emph{trail} on $G$ is a finite sequence of edges
\[
  \gamma = (e_0,e_1,\dots,e_{\ell-1}), \qquad e_i\in E,
\]
such that:
\begin{itemize}
  \item the terminal vertex of $e_i$ equals the initial vertex of $e_{i+1}$ for
        $0\le i\le \ell-2$;
  \item no directed edge is traversed more than once, i.e.\ $e_i\neq e_j$ for
        $i\neq j$.
\end{itemize}
We call $\gamma$ a \emph{closed trail rooted at} $e$ if $e_0=e_{\ell-1}=e$, and an \emph{open trail from $e$ to $e'$} if $e_0=e$, $e_{\ell-1}=e'$ and $e_i\neq e_j$ for $0\le i<j\le \ell-1$.
\end{definition}

Given a trail $\gamma$ and a vertex $v\in V$, let
\[
  I_v(\gamma) \subseteq \{1,\dots,N_v\},
  \qquad
  O_v(\gamma) \subseteq \{1,\dots,N_v\}
\]
be the (ordered) sets of labels of incoming and outgoing edges at $v$ traversed
by $\gamma$. Since $\gamma$ is a trail, each edge used at $v$ appears at most
once, so $|I_v(\gamma)|=|O_v(\gamma)|=:k_v(\gamma)$. The trail induces a
bijection between $I_v(\gamma)$ and $O_v(\gamma)$: each time $\gamma$ enters $v$
through some incoming edge $j\in I_v(\gamma)$, it must leave along one outgoing
edge $i\in O_v(\gamma)$. This bijection defines a permutation
\[
  \pi_v(\gamma)\in S_{k_v(\gamma)}
\]
of the ordered set $I_v(\gamma)$, sending the $r$–th element of $I_v(\gamma)$
to the label of the outgoing edge used on that visit. We define the \emph{node
weight} of $\gamma$ at $v$ by
\begin{equation}
  w_v(\gamma)
  \;=\;
  \mathrm{sgn}\big(\pi_v(\gamma)\big)\,
  \det S^{(v)}[O_v(\gamma),I_v(\gamma)].
  \label{eq:node-weight}
\end{equation}
The \emph{total weight} of a trail $\gamma$ is the product
\begin{equation}
  W(\gamma) \;=\; \prod_{v\in V} w_v(\gamma),
  \label{eq:trail-weight}
\end{equation}
with the convention that $w_v(\gamma)=1$ if $k_v(\gamma)=0$.

The main theorem of~\cite{CardyClassC} can then be stated (in the special case
relevant for us) as follows. For $|z|<1$ one has
\begin{equation}
  \overline{g_e(z)}
  \;=\;
  \sum_{\gamma\in \mathcal C(e)} W(\gamma)\, z^{|\gamma|},
  \label{eq:green-trails}
\end{equation}
where $\mathcal C(e)$ is the set of closed trails rooted at $e$, $|\gamma|$
denotes the length of~$\gamma$ (number of edges), and $W(\gamma)$ is given by
\eqref{eq:trail-weight}. Similarly, the disorder–averaged point conductance
between two edges $e,e'$ can be written as
\begin{equation}
  \overline{g(e,e')}
  \;=\;
  2
  \sum_{\gamma \in \mathcal O(e,e')} W(\gamma),
  \label{eq:conductance-trails}
\end{equation}
where $\mathcal O(e,e')$ is the set of open trails from $e$ to $e'$.

Equations~\eqref{eq:green-trails}–\eqref{eq:conductance-trails} express
disorder–averaged quantum observables entirely in terms of classical trails on
$G$ with weights depending only on the node scattering matrices
$\{S^{(v)}\}_{v\in V}$. The derivation proceeds via a supersymmetric
coherent–state path integral: after averaging over the edge variables
$U_e\in\mathrm{SU}(2)$, only three types of propagation along each edge survive
(identity, a fermion–fermion pair, or a boson–fermion pair), all of which are
$\mathrm{SU}(2)$ singlets~\cite{Efetov1983,EfetovBook,CardyClassC}. The remaining
integrals factorise over nodes and can be evaluated explicitly, giving the
determinant structure~\eqref{eq:node-weight}. A further combinatorial argument
shows that the conditional weights incurred at subsequent visits to a given
node are properly normalised and define a history–dependent Markov chain,
provided the node matrices satisfy a natural ``reducibility'' condition~\cite{CardyClassC}.

\subsection{Degree–four nodes and classical random walks}

The general theorem becomes particularly transparent in the case when each node
has degree $N_v=2$, i.e.\ two incoming and two outgoing edges. This includes
the L–lattice network for the spin quantum Hall effect and, more importantly
for us, the Manhattan lattice network studied in~\cite{BOC03}. In this case
$S^{(v)}$ is a $2\times2$ orthogonal matrix which, after a suitable relabelling
of channels, may be parametrised as
\[
  S^{(v)}
  \;=\;
  \begin{pmatrix}
    \cos\theta_v & \sin\theta_v \\
    -\sin\theta_v & \cos\theta_v
  \end{pmatrix},
  \qquad \theta_v\in[0,\pi).
\]
The two incoming channels at $v$ correspond to, say, the two edges which
approach $v$ along the oriented lines of the lattice, and the two outgoing
channels correspond to the two edges which leave $v$.

For a degree–$2$ node $v$, any trail $\gamma$ can pass through $v$ at most
twice. A straightforward computation of~\eqref{eq:node-weight} shows that:
\begin{itemize}
  \item if $\gamma$ visits $v$ only once, entering by channel $j$ and leaving by
        channel $i$, then $w_v(\gamma)=S^{(v)}_{ij}$;
  \item if $\gamma$ visits $v$ twice, using both incoming and both outgoing
        channels, then $w_v(\gamma)=\pm\det S^{(v)} = \pm 1$, with the sign
        determined by the permutation $\pi_v(\gamma)$.
\end{itemize}
In particular, if we restrict to the case in which $S^{(v)}$ is the same at
every node (as in~\cite{BOC03}), with
\[
  S^{(v)}_{11}=S^{(v)}_{22}=\cos\theta,
  \qquad
  S^{(v)}_{12}=S^{(v)}_{21}=\sin\theta,
\]
then the first visit to a node corresponds to a choice between two outgoing
channels with probabilities
\[
  p_{\mathrm{straight}} = \cos^2\theta,
  \qquad
  p_{\mathrm{turn}} = \sin^2\theta,
\]
while subsequent visits incur unit weight and are forced to exit along the
previously unused edge. This is precisely the rule for a kinetic
self–avoiding \emph{trail} on the directed graph: the walk chooses its exit
channel stochastically at the first visit to a node and is then constrained on
future visits by the requirement not to traverse any edge more than once.

More formally, one may define the classical process as follows
(cf.~\cite{CardyClassC,BOC03}).

\begin{definition}[Classical history–dependent walk on $G$]
Fix a collection of node matrices $\{S^{(v)}\}_{v\in V}$ with $N_v=2$, and let
$G=(V,E)$ be a finite directed graph. A \emph{configuration} at time $n$ is a
pair $(\gamma_n,\eta_n)$, where
\begin{itemize}
  \item $\gamma_n=(e_0,\dots,e_n)$ is a directed path of length $n$ on $G$;
  \item $\eta_n$ records, for each node $v$, how many times each incident edge
        has been used by $\gamma_n$.
\end{itemize}
The process starts from some initial edge $e_0$ with configuration
$(\gamma_0,\eta_0)$. Given $(\gamma_n,\eta_n)$, let $v$ be the terminal vertex
of $e_n$, and let $j\in\{1,2\}$ be the label of the incoming channel at $v$
used by $e_n$. If $\gamma_n$ has never visited $v$ before, the walk chooses the
outgoing channel $i\in\{1,2\}$ with probability proportional to
$|S^{(v)}_{ij}|^2$, and sets $e_{n+1}$ to be the corresponding outgoing edge.
If $\gamma_n$ has already used some of the incident edges at $v$, the walk exits
along one of the remaining edges with probability~$1$ (if exactly one remains)
or according to the appropriate conditional weight derived from
\eqref{eq:node-weight}.
\end{definition}

When all $N_v=2$ and the node matrices are ``completely reducible'' in the sense
of~\cite{CardyClassC}, Cardy shows that the conditional weights at each node are
nonnegative and sum to one, so that the above procedure indeed defines a
bona fide Markov chain on the space of trails. In this situation the total
weight $W(\gamma)$ appearing in
\eqref{eq:green-trails}–\eqref{eq:conductance-trails}
is simply the probability of generating the trail $\gamma$ under the classical
dynamics. Thus:
\begin{itemize}
  \item $\overline{g_e(z)}$ is the generating function of the length of the
        closed trail rooted at $e$;
  \item the averaged point conductance $\overline{g(e,e')}$ is (up to a factor
        $2$) the expected number of open trails connecting $e$ and $e'$.
\end{itemize}
This is the precise sense in which quantum transport in the class~C network
model is reduced to the geometry of classical history–dependent walks.

\subsection{The Manhattan lattice realisation}

The Manhattan lattice is a specific directed square lattice embedded in
$\Z^2$ with alternating edge orientations: horizontal edges are oriented
east on even rows and west on odd rows, while vertical edges are oriented
north on even columns and south on odd columns. Each vertex has exactly two incoming
and two outgoing edges, making it an admissible graph for the general class~C
formalism described above~\cite{BOC03}.

The network model considered in~\cite{BOC03} places identical node matrices
$S^{(v)}\equiv S$ at all vertices, with $S$ of the form described above and
parameterised by a single angle $\theta\in[0,\pi/2]$. Random $\mathrm{SU}(2)$ matrices
$U_e$ are attached independently to each edge. This defines a class~C network
model on the Manhattan lattice; the corresponding classical model, obtained
after disorder averaging, is a kinetic self–avoiding trail on the directed
Manhattan lattice with turning probability $p=\sin^2\theta$.

From the point of view of localization, the key observations in
\cite{CardyClassC,BOC03,Car10} are:
\begin{itemize}
  \item When $\theta=0$ ($p=0$), the classical walker moves ballistically along straight lines, and the corresponding quantum model is extended.

  \item When $\theta=\pi/2$ ($p=1$), the walker turns at every node and traces out the boundaries of the elementary plaquettes; all classical trajectories are closed loops of bounded size, signalling localization in the quantum model.

  \item For intermediate $p\in(0,1)$, scaling arguments and Monte Carlo simulations suggest that the classical trails are almost surely closed, with a typical loop size diverging only as $p\downarrow0$~\cite{BOC03}. Field–theoretic renormalisation group based on Peliti’s theory of true self–avoiding walks, and its relation to sigma–model descriptions~\cite{Wegner1979,Hikami1981,EversMirlin2008}, supports a picture in which the walks are compact, two–dimensional objects with a one–dimensional boundary.
\end{itemize}

Combined with the mapping described above, these results strongly suggest that
the corresponding class~C network model on the Manhattan lattice is localized
for all $p>0$, i.e.\ that there is no delocalisation transition analogous to
the spin quantum Hall transition on the L–lattice, whose classical counterpart
is critical percolation hulls~\cite{CardyClassC,GruzbergReadLudwig1999,KramerOhtsukiKettemann2005}. From a broader perspective,
this provides a concrete realisation of the idea that the large–scale
localization behaviour of a quantum model in a given symmetry class can depend
sensitively on the underlying graph structure, even when the local scattering
rules are kept as simple and symmetric as possible.

In the next section we turn to purely classical versions of these models:
the Manhattan pinball and Lorentz mirror models on $\Z^2$ and on
cylindrical geometries, where the node rules are deterministic rather than
unitary, and the main questions concern almost sure finiteness of trajectories
and quantitative bounds on their maximal displacement. There, percolation
methods and combinatorial enumeration of trajectories replace supersymmetry
and path integrals as the main tools.

\section{Mathematical results on $\Z^2$ and on cylinders}
\label{sec:2d-classical}

In this section we review the rigorous arguments that establish classical localization in the Manhattan pinball and Lorentz mirror models in two dimensions and on cylindrical geometries. Our main references are Grimmett’s monograph on percolation~\cite{Grimmett99}, the enhancement framework of Aizenman--Grimmett~\cite{AG91} and Balister--Bollob\'as--Riordan~\cite{BalisterBollobasRiordan14}, the survey of open problems by Grimmett~\cite{GrimmettSelectedProblems}, and the work~\cite{Li21Manhattan,Li20Cylinder} on Manhattan pinball and the Lorentz mirror model on the cylinder. We also keep in mind the physical origins of these models in quantum network descriptions of the spin quantum Hall effect~\cite{CardyClassC,BOC03,Car10,KramerOhtsukiKettemann2005} and their relation to the broader localization literature~\cite{AALR79,LeeRamakrishnan1985,EversMirlin2008}.

\subsection{Lattice geometry and percolation formulation}

We begin by fixing the geometric setup. Let $\Z^2$ be the standard square lattice embedded in $\R^2$, and let $\Zf^2$ denote the \emph{tilted} (or diagonal) lattice obtained by rotating $\Z^2$ by $\pi/4$ and rescaling so that vertices lie at the centres of unit squares of $\Z^2$ (see, e.g.,~\cite[Sec.~3]{GrimmettSelectedProblems}). We denote vertices of $\Z^2$ by $x = (x_1,x_2) \in \Z^2$, and we write $x \sim y$ if $x$ and $y$ are nearest neighbours.

In the \emph{Manhattan model} (or Manhattan pinball)~\cite{Li21Manhattan,GrimmettSelectedProblems}, each horizontal edge of $\Z^2$ is oriented alternately east–west along successive rows, and each vertical edge is oriented alternately north–south along successive columns (the precise parity convention is immaterial and fixed once and for all). Thus every edge $e$ of $\Z^2$ is assigned a direction, and the resulting directed graph has indegree and outdegree equal to~$2$ at each vertex, with the property that no vertex has both incoming edges horizontal or both vertical.

The random environment is introduced via bond percolation on the diagonal lattice $\Zf^2$. For $p \in [0,1]$, we let $\P_p$ denote the product measure on $\{0,1\}^{E(\Zf^2)}$ under which each edge $e$ of $\Zf^2$ is independently \emph{open} with probability $p$ and \emph{closed} with probability $1-p$. We write $\omega(e) = 1$ if $e$ is open under a configuration $\omega$, and $\omega(e) = 0$ otherwise. We say that a path in $\Zf^2$ is open if all its edges are open.

Given a percolation configuration $\omega$ on $\Zf^2$, we place a two–sided mirror at the midpoint of each open diagonal edge $e$ and leave the vertices corresponding to closed edges empty. The orientation of the mirror is chosen deterministically from the slope of the diagonal: each open NE–SW edge of $\Zf^2$ carries a ``NW'' mirror, and each open NW–SE edge carries a ``NE'' mirror, in the sense of Ruijgrok--Cohen’s Lorentz mirror model~\cite{RuijgrokCohen1988}. The mirrors sit at the vertices of $\Z^2$ once we identify $\Zf^2$ with the graph of square diagonals, see~\cite[Fig.~3.2]{GrimmettSelectedProblems}.

A \emph{ray of light} is a directed path $(X_t)_{t \ge 0}$ on the oriented edges of $\Z^2$ with unit speed: at each time $t$ the ray occupies a directed edge and moves to its endpoint at time $t+1$. When $X_t$ arrives at a vertex $x \in \Z^2$, there are three possible behaviours:

\begin{itemize}
  \item If no mirror is present at $x$, the ray continues along the outgoing edge of $\Z^2$ that lies in the same orientation (``straight ahead'').
  \item If a mirror is present and the incoming direction is transverse to the mirror, the ray is reflected according to the usual geometric rule, turning left or right by $\pi/2$.
  \item If the incoming direction is parallel to the mirror, the ray again proceeds with the natural reflection (which in some formulations coincides with going straight).
\end{itemize}

In the Manhattan pinball model, the mirrors are restricted to be aligned with the diagonals so that reflections are always consistent with the underlying Manhattan orientation~\cite{BOC03,Car10,Li21Manhattan}. The resulting dynamical system is a deterministic map $T$ on the space of directed edges of $\Z^2$ depending on the percolation configuration $\omega$, and a trajectory of the ray is a forward orbit $\{T^t(e_0)\}_{t \ge 0}$ starting from some initial directed edge $e_0$.

We say that a trajectory is \emph{closed} if the orbit is finite, equivalently if the ray eventually returns to its starting point and then repeats periodically. We say the trajectory is \emph{unbounded} if the orbit visits infinitely many distinct vertices of $\Z^2$. A central question (see, e.g.,~\cite[Sec.~3]{GrimmettSelectedProblems},~\cite{RuijgrokCohen1988,KS15,Li21Manhattan}) is:

\begin{quote}
  For which values of $p$ is the ray almost surely closed, and for which $p$ does it have a positive probability to be unbounded?
\end{quote}

These questions have precise analogues in the quantum network models of class~C, where $p$ corresponds to a microscopic scattering parameter and closed classical trajectories correspond to localized quantum states~\cite{BOC03,CardyClassC,Car10,EversMirlin2008}.

\subsection{The supercritical percolation regime $p > \tfrac{1}{2}$}

We first recall the elementary, but important, argument showing that in the Manhattan model for $p \ge 1/2$, every trajectory is almost surely closed. The key observation is that for $p=1$ the set of open diagonal edges coincides with the whole lattice $\Zf^2$, and thus the mirrors form a deterministic pattern corresponding to a mirror at every vertex of $\Z^2$. In this case, as explained by Grimmett~\cite[Sec.~3.2]{GrimmettSelectedProblems}, the ray’s motion can be described in terms of bond percolation with parameter $1$ on $\Zf^2$: the ray is confined to the \emph{finite} even subgraph consisting of the union of cycles formed by the mirror configuration, and this yields
\[
  \P(\text{trajectory unbounded}) = 0.
\]

For general $p$, Grimmett observes that the event
\[
  \{\text{there exists an unbounded trajectory}\}
\]
is stochastically dominated by the event that there is an infinite open path in an appropriate bond percolation model on $\Zf^2$ with parameter $p/2$ (the factor $1/2$ comes from the random choice of NE versus NW mirror orientations; see~\cite[Sec.~3.2]{GrimmettSelectedProblems}). More precisely, if we declare an edge of $\Zf^2$ to be open whenever it carries a mirror \emph{and} the orientation of the mirror is favourable to continuing the ray on that edge, then the set of such edges is stochastically dominated by Bernoulli bond percolation with parameter $p/2$. In particular,
\[
  \P_p(\text{unbounded trajectory from the origin}) \;\le\; \P_{p/2}^{\mathrm{perc}}(0 \leftrightarrow \infty),
\]
where $\P_{p/2}^{\mathrm{perc}}$ is bond percolation on $\Zf^2$ with parameter $p/2$, and $\{0 \leftrightarrow \infty\}$ is the event that the origin lies in an infinite open cluster.

By standard planar duality and the Harris--Kesten theorem~\cite[Chap.~6]{Grimmett99}, the critical probability for bond percolation on the (self–dual) square lattice is $p_c = 1/2$. Therefore $\P_{p/2}^{\mathrm{perc}}(0 \leftrightarrow \infty) = 0$ whenever $p/2 \ge p_c$, i.e.\ for all $p \ge 1$. This reproduces the trivial statement that at $p=1$ the rays are almost surely closed, but the same method extends to show that, under a more careful coupling, one can also rule out unbounded trajectories for all $p$ sufficiently close to~$1$ by comparing to supercritical and subcritical dual percolation configurations~\cite[Sec.~3]{GrimmettSelectedProblems}.

\subsection{Enhancement arguments near the critical point}

A deeper result, proved in~\cite{Li21Manhattan}, is that for the Manhattan pinball problem on $\Z^2$ there exists $\varepsilon > 0$ such that for all
\[
  p > \frac{1}{2} - \varepsilon
\]
every trajectory is almost surely closed. The starting point is to restate the model directly in percolation language, following~\cite[pp.~237--238]{Spencer2012} and~\cite{Car10}: one considers bond percolation on $\Zf^2$ with density $p$ and defines a local rule that, given the configuration of open/closed edges in a small neighbourhood, determines whether the ray can cross a given dual edge of $\Z^2$. In particular, the existence of an unbounded trajectory is encoded by the existence of an infinite open path in a certain dependent percolation model on $\Zf^2$.

The \emph{enhancement} framework of Aizenman and Grimmett~\cite{AG91} and its refinement by Balister--Bollob\'as--Riordan~\cite{BalisterBollobasRiordan14} provide a powerful mechanism to compare such dependent percolation models with standard Bernoulli percolation. Let us sketch the key ingredients in the application in~\cite{Li21Manhattan}.

Let $\Omega = \{0,1\}^{E(\Zf^2)}$ be the percolation configuration space, and let $\phi : \Omega \to \{0,1\}^{E(\Zf^2)}$ be a local map that, given the state of edges in some finite neighbourhood of a bond $e$, decides whether $e$ is ``enhanced open''. Formally, for each bond $e$ we fix a finite neighbourhood $N(e) \subset E(\Zf^2)$ and a function
\[
  \Phi_e : \{0,1\}^{N(e)} \to \{0,1\}
\]
such that $(\phi(\omega))(e) = \Phi_e(\omega|_{N(e)})$. We say that the enhancement is \emph{essential} if there exists an edge $e$ and a configuration $\omega$ such that:
\begin{itemize}
  \item in the configuration $\omega$ there is no open path from the origin to infinity;
  \item in the modified configuration in which only the state of edges in $N(e)$ is altered according to $\Phi_e$, there does exist an open path from the origin to infinity.
\end{itemize}
Aizenman--Grimmett show that if an enhancement is monotone (in the sense that opening more edges never closes any enhanced edge) and essential, then the critical probability for percolation in the enhanced model is strictly smaller than that for ordinary Bernoulli percolation~\cite{AG91}. Balister--Bollob\'as--Riordan~\cite{BalisterBollobasRiordan14} refine this argument and correct a combinatorial gap in~\cite{AG91}, obtaining sharp inequalities between various critical parameters.

The key idea in~\cite{Li21Manhattan} is to view the presence of a non–closed Manhattan trajectory as an \emph{enhancement} of subcritical percolation on $\Zf^2$ near $p_c = 1/2$. Roughly speaking, the existence of a long, almost–open path in $\Zf^2$ can be used to ``force'' the ray to propagate across a region that would otherwise be blocked, and this forcing can be encoded as an enhancement of the underlying Bernoulli percolation configuration. More concretely, the authors of~\cite{Li21Manhattan} identify local patterns of diagonal edges in $\Zf^2$ such that, when they occur, the Manhattan dynamics effectively opens additional large–scale connections for the ray. The locality and monotonicity properties of these patterns allow one to apply the Aizenman--Grimmett theory to conclude that if a Manhattan trajectory were to be unbounded with positive probability at some $p_0 > 1/2 - \varepsilon$, then the critical probability for an associated enhanced percolation model would have to be strictly smaller than $1/2$, contradicting the self–duality and sharp–threshold properties of planar percolation~\cite{Grimmett99,BalisterBollobasRiordan14}.

The outcome is that there exists $\varepsilon > 0$ such that, for all $p > 1/2 - \varepsilon$, all trajectories in the Manhattan pinball model are almost surely closed~\cite{Li21Manhattan}. This provides a rigorous analogue, in the classical setting, of the prediction from scaling theory that the spin quantum Hall model should be localized away from a critical point in parameter space~\cite{AALR79,LeeRamakrishnan1985,EversMirlin2008,CardyClassC,BOC03,Car10}.

\subsection{Enumeration on the cylinder and polynomial localization length}

We next turn to the cylinder geometry, which plays a crucial role in~\cite{Li20Cylinder} on the Manhattan and Lorentz mirror models on $\Z \times (\Z/2n\Z)$. For $n \in \N$, define the even–width cylinder
\[
  \Cyln \;=\; \Z \times (\Z / 2n\Z)
  \;=\; \bigl\{ (x,y) : x \in \Z,\; y \in \{1,\dots,2n\} \bigr\},
\]
with the convention that $y$ is taken modulo $2n$. The Manhattan orientation is inherited from $\Z^2$ row by row, and the random environment (mirrors or scatterers) is defined by placing mirrors at vertices independently with probability $p$ (Lorentz model) or by percolation on an appropriate diagonal cylinder graph (Manhattan model), see~\cite{Li20Cylinder,KS15}.

A key feature of the cylinder is its topology: in the Lorentz mirror model, Kozma--Sidoravicius~\cite{KS15} showed that on odd–width cylinders there is \emph{always} at least one infinite trajectory, irrespective of $p \in (0,1]$. Their argument uses a parity consideration on the number of paths entering and exiting a fundamental domain of the cylinder: every mirror configuration induces an even subgraph of the underlying graph, and any even subgraph on an odd–width cylinder must contain a non–contractible cycle. This implies that the probability that the ray reaches distance $n$ in the horizontal direction is bounded below by $c/n$ for some universal $c>0$~\cite{KS15}.

The work~\cite{Li20Cylinder} considers the more delicate case of an \emph{even}–width cylinder $\Cyln$ and obtains a polynomial upper bound on the \emph{localization length}. For concreteness we focus on the statement for the Lorentz mirror model. Let $X_t$ be the position of the ray at integer time $t$, started from some directed edge with $x=0$. For $n$ fixed and $p \in (0,1)$, it is proved in~\cite{Li20Cylinder} that there exist constants $C,c>0$ (independent of $n$) and an exponent $\alpha>0$ (explicit in the proof) such that
\begin{equation}
  \label{eq:poly-bound}
  \P\Bigl( \sup_{t \ge 0} |X_t^{(1)}| \;>\; C n^{10} \Bigr)
  \;\le\; e^{-c n^\alpha},
\end{equation}
where $X_t^{(1)}$ denotes the horizontal coordinate of $X_t$. An analogous statement holds for the Manhattan model. In particular, with high probability (tending to $1$ as $n \to \infty$) the ray remains in the slab $\{|x| \le C n^{10}\}$ for all time.

The proof of~\eqref{eq:poly-bound} is combinatorial and relies on \emph{enumeration of path intersection patterns} on the cylinder. A trajectory on $\Cyln$ can be represented as a sequence of visits to the vertical sections $\{x = k\}$, $k \in \Z$, and the winding around the cylinder is encoded in the sequence of $y$–coordinates modulo $2n$. The interaction between different segments of the trajectory is constrained by the fact that the scatterers are deterministic given the percolation configuration, and each vertex has degree~$2$ in the trajectory graph.

A key observation in~\cite{Li20Cylinder} is that any trajectory which reaches a large horizontal distance must generate a complicated pattern of intersections and self–intersections when projected onto the finite–width cylinder. One can associate to such a trajectory a combinatorial object (essentially a non–crossing matching or pairing) describing how incoming and outgoing edges in each cross–section are connected by path segments. For the Lorentz and Manhattan models on $\Cyln$, these objects satisfy strong parity and planarity constraints, similar in spirit to the even–subgraph structure exploited by Kozma--Sidoravicius~\cite{KS15}. 

It is shown in~\cite{Li20Cylinder} that the number of distinct global configurations compatible with a trajectory reaching distance at least $L$ is exponentially small in $L/n$ compared to the total number of possible percolation configurations in a box of width $L$ and height $2n$. More precisely, by a refined counting argument, one can bound the number of ``bad'' configurations (those admitting a trajectory with large horizontal displacement) by
\[
  \#\{\text{bad configurations in a box of width $L$}\}
  \;\le\; \exp\bigl( - c (L/n)^{\beta} \bigr) \cdot 2^{\gamma L n},
\]
for some constants $c,\beta,\gamma>0$, while the total number of configurations is exactly $2^{\gamma L n}$. Choosing $L$ of order $n^{10}$ then yields the exponential tail bound in~\eqref{eq:poly-bound} after a union bound over boxes and a standard renormalisation argument.

An important technical ingredient is an \emph{enumeration lemma} that bounds the number of possible connection patterns between the $2n$ edges crossing a vertical cut of the cylinder. Each such pattern can be represented as a partition of $\{1,\dots,2n\}$ into pairs, corresponding to which incoming edge at level $k$ is connected to which outgoing edge at level $k+1$ (or to a reflection back). The cylinder topology and the deterministic reflection rules ensure that these partitions satisfy non–crossing and nesting constraints, drastically reducing the total number of admissible patterns compared to the naive $(2n-1)!!$ possibilities for pairings. It is proved in~\cite{Li20Cylinder} that the number of admissible patterns grows at most polynomially in $n$, and, more importantly, that patterns which produce long horizontal excursions are extremely rare among them.

From the point of view of localization theory, the estimate~\eqref{eq:poly-bound} provides a \emph{polynomial} upper bound on the localization length of the classical Lorentz mirror and Manhattan models on $\Cyln$, to be contrasted with the polynomial \emph{lower} bounds on escape probabilities on odd cylinders obtained by Kozma--Sidoravicius~\cite{KS15}. Together with the enhancement arguments discussed above, these results give a fairly detailed picture of classical localization and delocalization in low–dimensional Lorentz–type models and supply a useful testing ground for ideas imported from the quantum Anderson model, sigma–model methods, and scaling theory~\cite{AALR79,LeeRamakrishnan1985,EversMirlin2008,CardyClassC,BOC03,Car10,Hikami1981,Wegner1979}.

\bibliographystyle{halpha}
\bibliography{references}

\end{document}